\documentclass{amsart}

\newtheorem{theorem}{Th\'eor\`eme}

\numberwithin{theorem}{section} \theoremstyle{plain}

\newtheorem{conjecture}[theorem]{Conjecture}

\newtheorem{corollary}[theorem]{Corollaire}

\newtheorem{lemma}[theorem]{Lemme}

\newtheorem{proposition}[theorem]{Proposition}

\theoremstyle{definition}

\newtheorem{definition}[theorem]{D\'efinition}

\newtheorem*{notation}{Notation}

\newtheorem{example}[theorem]{Exemple}

\newtheorem{remark}[theorem]{Remarque}

\numberwithin{equation}{section}

\DeclareMathOperator{\pic}{Pic}

\newcommand{\xx}{\mathcal{X}}

\newcommand{\pp}{\mathfrak{p}}

\newcommand{\kappap}{\kappa_{\pp}}

\newcommand{\qp}{q_{\pp}}

\newcommand{\qq}{\mathbb{Q}}

\newcommand{\zz}{\mathbb{Z}}

\newcommand{\qqell}{\qq_{\ell}}

\DeclareMathOperator{\rg}{rang} \DeclareMathOperator{\tr}{Tr}

\newcommand{\ssc}{\mathcal{S}}

\DeclareMathOperator*{\res}{R\acute{e}s}
\DeclareMathOperator*{\ord}{ord} \DeclareMathOperator{\ns}{NS}

\DeclareMathOperator{\ddiv}{Div}

\newcommand{\ppb}{\mathbb{P}}

\newcommand{\aab}{\mathbb{A}}

\DeclareMathOperator{\alb}{Alb} \DeclareMathOperator{\kker}{Ker}

\DeclareMathOperator{\gal}{Gal} \DeclareMathOperator{\spec}{Spec}
\DeclareMathOperator{\et}{\acute{e}t}

\newcommand{\gk}{G_k}

\newcommand{\xbar}{\overline{\xx}}

\newcommand{\ov}{\overline}\newcommand\jour{\number\day\space
\space\ifcase\month\or janvier\or f\'evrier\or

mars\or avril\or mai\or juin\or

juillet\or ao\^ut\or septembre\or octubre\or novembre\or
d\'ecembre\fi\space \space\number\year}

\newcommand{\unc}{\underline{c}}\newcommand{\cac}{\mathcal{C}}

\begin{document}

\title{Sur le rang des Jacobiennes sur un corps de fonctions}

\author{Marc Hindry}

\address{Universit\'e Denis Diderot Paris VII\\U.F.R. Math\'ematiques\\
   case 7012\\ 2 Place Jussieu\\ 75251 Paris,
France} \email{hindry@math.jussieu.fr}

\author{Am\'{\i}lcar Pacheco}

\address{Universidade Federal do Rio de Janeiro (Universidade do
Brasil)\\ Departamento de Ma\-te\-m\'a\-ti\-ca Pura\\ Rua Guai\-aquil
83, Cachambi, 20785-050 Rio de Janeiro,
RJ, Brasil} \email{amilcar@impa.br}\thanks{\textbf{Remerciements.}
Am\'{\i}lcar Pacheco a \'et\'e
partiellement soutenu par
la bourse de recherche CNPq 300896/91-3 et par le projet PRONEX
41.96.0830.00. Ce travail a commenc\'e pendant
une visite du premier auteur au deuxi\`eme auteur \`a l'IMPA dans le
cadre de l'accord
Br\'esil-France  69.0014/01-5, les
deux auteurs remercient cet accord pour son soutien financier et
aussi l'IMPA pour sa chaleureuse
ambiance scientifique. Les auteurs remercient aussi le rapporteur
pour ses suggestions, notamment
concernant la preuve du Lemme \ref{geom1}.}

\date{\jour}

\begin{abstract}Soit $f:\mathcal{X}\rightarrow C$ une surface
projective fibr\'ee au dessus d'une courbe et d\'efinie sur
   un corps de nombres $k$. Nous donnons une interpr\'etation du rang
du groupe de Mordell-Weil sur $k(C)$ de la  jacobienne de la fibre
g\'en\'erique (modulo la partie constante) en
termes de moyenne des traces de Frobenius sur les fibres de $f$.
L'\'enonc\'e fournit une
r\'einterpr\'etation de la
conjecture de Tate pour la surface $\mathcal{X}$  et g\'en\'eralise
des r\'esultats de Nagao, Rosen-Silverman et Wazir.

\par\bigskip

\noindent\textsc{Abstract.} \textbf{On the rank of Jacobians over
function fields.} Let $f:\mathcal{X}\rightarrow C$
be a projective surface fibered over a curve and defined over a
number field $k$. We give an
interpretation of the rank of
the  Mordell-Weil group over $k(C)$ of the jacobian of the generic
fibre  (modulo the constant part) in terms
of average of the traces of Frobenius on the fibers of $f$. The
results also give a reinterpretation
of the
   Tate conjecture for the surface $\mathcal{X}$  and generalizes
results of Nagao, Rosen-Silverman and Wazir.\end{abstract}

   \maketitle

\section{Introduction}

Soient $k$ un corps de nombres, $\mathcal{X}$ une surface projective
lisse irr\'eductible sur $k$, $C$ une
courbe projective lisse irr\'eductible sur $k$, $f:\mathcal{X}\to C$
un morphisme propre plat tel
que les fibres soient des
courbes de genre arithm\'etique $g\ge 1$; ces hypoth\`eses
entra\^{\i}nent que la fibre g\'en\'erique est lisse
et irr\'eductible. Soient $K=k(C)$ le corps de fonctions de $C$,
$X/K$ la fibre g\'en\'erique de
$f$, $J_X$ la vari\'et\'e
Jacobienne de $X$ et $(\tau,B)$ la ${K}/{k}$-trace de $J_X$. Un
th\'eor\`eme de Lang-N\'eron affirme
que
$J_X(k(C))/\tau(B(k))$ et m\^eme $J_X(\bar{k}(C))/\tau(B(\bar{k}))$
sont des groupes de type fini. Shioda \cite{sh1}
a donn\'e une interpr\'etation du deuxi\`eme groupe comme quotient du
groupe de N\'eron-Severi de la
surface
$\mathcal{X}$ par un sous-groupe explicite.  Nous allons discuter une
interpr\'etation du {\it rang} du premier
groupe en termes de l'arithm\'etique de la surface $\mathcal{X}$.

Soit $S$ un ensemble fini d'id\'eaux premiers de l'anneau d'entiers
$\mathcal{O}_k$ de $k$ tel que pour
tout
$\mathfrak{p}\notin S$, $\mathcal{X}$ et $C$  aient bonne
r\'eduction  modulo $\pp$ et que la
r\'eduction
$f_{\mathfrak{p}}$ de $f$ modulo $\mathfrak{p}$ soit un morphisme
propre et plat
$f_{\mathfrak{p}}:\mathcal{X}_{\mathfrak{p}}\to C_{\mathfrak{p}}$
ayant pour fibres des courbes de genre
arithm\'etique
$g$ sur le corps r\'esiduel $\kappa_{\mathfrak{p}}$ de $\mathfrak{p}$
avec $q_{\mathfrak{p}}$ \'el\'ements. Pour
chaque
$c\in C_{\pp}$, soit $\xx_{\pp,c}=f_{\pp}^{-1}(c)$ la fibre de
$f_{\pp}$ en $c$.

\begin{notation}Soient $l$ un corps, $\overline{l}$ une cl\^oture
alg\'ebri\-que de $l$ et $V$ une vari\'et\'e alg\'ebri\-que projective lisse d\'efinie sur $l$. Notons
$\overline{V}=V\times_l\overline{l}$. Si $l$ est parfait, notons aussi $G_l=\gal(\overline{l}/l)$.

\end{notation}

Soit $F_{\pp}\in G_k$ un \'el\'ement de Frobenius et $I_{\pp}\subset
G_k$ son groupe d'inertie. Soit
$\Delta=\{c\in C\,|\,\xx_{\pp,c}$ est singuli\`ere$\}$ le lieu
discriminant de $f$. Apr\`es avoir
\'elargi $S$, si n\'ecessaire, on
peut supposer que pour tout $\pp\notin S$ le lieu discriminant
$\Delta_{\pp}$ de $f_{\pp}$ soit \'egal \`a
la r\'eduction de $f$ modulo $\pp$.

Soit $\ov{F}_{\pp}$ l'automorphisme de Frobenius sur $H^1_{\et}(\ov{\xx}_{\pp,c},\qq_{\ell})$ induit par
l'automor\-phis\-me de Frobenius g\'eom\'etrique de $\kappa_{\pp}$. Si $c\in (C_{\pp}-\Delta_{\pp})(\kappa_{\pp})$, on
d\'efinit $a_{\pp}(\xx_{\pp,c})=\tr(\ov{F}_{\pp}\,|\,H^1_{\et}(\ov{\xx}_{\pp,c}, \qq_{\ell})$. Si
$c\in\Delta_{\pp}(\kappa_{\pp})$, on remplace $H^1_{\et}$ par le groupe $H^1_{\unc}$ de cohomologie $\ell$-adique \`a
support propre, i.e., $a_{\pp}(\xx_{\pp,c})=\tr(\ov{F}_{\pp}\,|\,H^1_{\unc}(\ov{\xx}_{\pp,c} ,\qq_{\ell}))$.

Soit
$a_{\pp}(B)=\tr(F_{\pp}\,|\,H^1_{\et}(\ov{B},\qq_{\ell})^{I_{\pp}})$
. Apr\`es avoir ajout\'e un nombre
fini d'id\'eaux premiers de $k$ \`a $S$ on peut supposer que, pour
tout $\pp\notin S$, la
vari\'et\'e $B$ ait bonne
r\'eduction $B_{\pp}$ modulo $\pp$. Dans ce cas-l\`a,
$a_{\pp}(B)=\tr(\ov{F}_{\pp}\,|\,H^1_{\et}(\ov{B}_{\pp},\qq_{\ell}))$.
On d\'efinit alors
$$A_{\mathfrak{p}}(\mathcal{X})=\frac 1{q_{\mathfrak{p}}}\sum_{c\in
C_{\pp}(\kappa_{\mathfrak{p}})}a_{\mathfrak{p}}(\mathcal{X}_{\pp,c})$$
la {\it trace moyenne} de Frobenius ainsi
que
$$A_{\mathfrak{p}}^*(\mathcal{X})=A_{\mathfrak{p}}(\mathcal{X})-a_{\mathfrak {p}}(B),$$ que nous
appellerons {\it trace moyenne r\'eduite}.

Soit $L_2(\mathcal{X}/k,s)=\prod_{\mathfrak{p}}\det(1-F_{\mathfrak{p}}q_{\mathfrak{p}}^{-s}\,|\,
H^2_{\text{\'et}}(\xbar,\mathbb{Q}_{\ell})^{I_{\pp}})^{-1}$ la fonction $L$ de $\mathcal{X}$ associ\'ee \`a
$H^2_{\text{\'et}}(\xbar,\qq_{\ell})$. Soient $\pic(\mathcal{X})$ le groupe des classes de diviseurs de $\mathcal{X}$,
$\pic^0(\mathcal{X})$ le sous-groupe de diviseurs alg\'ebriquement \'equivalents \`a z\'ero,
$\ns(\mathcal{X})={\pic(\mathcal{X})}/{\pic^0(\mathcal{X})}$ le groupe de N\'eron-Severi de $\mathcal{X}$ et
$\ns(\xx/k)$ le sous-groupe de classes de diviseurs de $\ns(\xx)$ qui sont d\'efinies sur $k$. Notons que si $z$ est
fix\'e par $G_k$, i.e., $z\in\ns(\xx)^{G_k}$, il existe un multiple $nz$ de $z$ tel que $nz\in\ns(\xx/k)$. Donc,
$\ns(\xx)^{G_k}\otimes\qq\cong\ns(\xx/k)\otimes\qq$ et en particulier $\rg(\ns(\xx)^{G_k})=\rg(\ns(\xx/k))$.

\begin{conjecture}[Conjecture de Tate \cite{ta2}]\label{conjT}  La
fonction $L_2(\mathcal{X}/k,s)$  pos\-s\`e\-de \linebreak un p\^ole
en $s=2$ d'ordre
$\rg(\ns(\xx/k))$.
\end{conjecture}

\begin{remark}Il s'agit d'une version de la conjecture de Tate pour
les diviseurs, la conjecture g\'en\'erale concerne tous les cycles
alg\'ebriques.

\begin{enumerate}
\item On peut se dispenser de l'hypoth\`ese d'un prolongement
m\'eromorphe au voisinage de $s=2$ en interpr\'etant
la phrase ``$L_2(\mathcal{X}/k,s)$ poss\`ede un p\^ole d'ordre $t$ en
$s=2$" comme signifiant
$$\lim_{\Re(s)>2, s\to 2}(s-2)^tL_2(\mathcal{X}/k,s)=\alpha\not= 0.$$
De m\^eme, si $f(s)$ est holomorphe sur $\Re(s)>\lambda$ et si
$\lim_{s\to\lambda}(s-\lambda)f(s)=\alpha\not= 0$,
on appellera $\alpha$ le {\it r\'esidu} de la fonction $f(s)$ en
$s=\lambda$ et on \'ecrira
$\res_{s=\lambda}f(s)=\alpha$.
\item Dans la plupart des cas o\`u l'on sait d\'emontrer un
prolongement analytique de $L_2(\mathcal{X}/k,s)$ \`a
la droite $\Re(s)=2$, on sait \'egalement d\'emontrer que la fonction
ne s'annule pas sur cette
droite. Cette
propri\'et\'e est importante car elle permet d'appliquer un
th\'eor\`eme Taub\'erien \cite[Chapter
XV]{la} \`a la d\'eriv\'ee logarithmique de $L_2(\mathcal{X}/k,s)$.
\end{enumerate}

\end{remark}

Le but de ce papier est de prouver le th\'eor\`eme suivant.

\begin{theorem}\label{thmA}La Conjecture \ref{conjT} pour la surface
$\mathcal{X}$ implique
\begin{equation}\label{NagaoAn}\res_{s=1}\left(\sum_{\mathfrak{p}\notin
S}-A_{\mathfrak{p}}^*(\mathcal{X})\frac{\log(q_{\mathfrak{p}})}
{q_{\mathfrak{p}}^s}\right)=\rg\left(\frac{J_X(K)}{\tau
B(k)}\right).\end{equation} Sous
l'hypoth\`ese additionelle que $L_2(\mathcal{X}/k,s)$ se prolonge
analytiquement  sur
la droite $\Re(s)=2$ et n'a pas de z\'eros sur cette droite, on
conclut aussi que
\begin{equation}\label{NagaoAr}\lim_{T\to\infty}\frac
1T\left(\sum_{\substack{\mathfrak{p}\notin S\\ q_{\wp}\le
T}}-A_{\mathfrak{p}}^*(\mathcal{X})\log(q_{\mathfrak{p}})\right)=\rg\left(\frac{J_X(K)}{\tau
B(k)}\right).\end{equation}
\end{theorem}

\begin{remark} En fait l'\'egalit\'e (\ref{NagaoAn}) est
essentiellement \'equivalente \`a la conjecture de Tate (\ref{conjT}), voir le
corollaire~(\ref{corTa-Na}). L'\'egalit\'e (\ref{NagaoAn}) est une
g\'en\'eralisation de la
\emph{conjecture analytique de Nagao} \cite[Nagao's
Conjecture 1.1$^\prime$]{rosi} et (\ref{NagaoAr}) de la
\emph{conjecture Taub\'erienne de Nagao}
\cite[Nagao's Conjecture 1.1]{rosi} telles qu'elles sont formul\'ees
par Rosen et Silverman pour
les fibrations de genre 1 admettant
une section; lorsque de plus $C=\ppb^1$ et $k=\qq$, la forme
Taub\'erienne  est la conjecture originelle, due \`a
Nagao
\cite{nag1,nag2} (nous pr\'ef\'erons appeler ``Taub\'erienne" la
forme de la conjecture que
Rosen et Silverman appellent ``{\it arithmetic}").\end{remark}

\begin{remark}Le premier travail dans la direction du Th\'eor\`eme
\ref{thmA} est d\^u \`a Nagao \cite{nag2} qui a formul\'e  la
conjecture Taub\'erienne \`a partir de
calculs explicites dans le cas o\`u $g=1$ et $f$ admet une section
(donc $X$ est une courbe
elliptique)  et $K=\qq(x)$. Son
objectif, atteint dans plusieurs cas, \'etait de produire des courbes
elliptiques sur $\qq$ avec rang
``\emph{assez grand}''. Par la suite, Rosen et Silverman dans
\cite[Theorem 1.3]{rosi} on trait\'e
le cas g\'en\'eral o\`u $K=k(C)$
est le corps de fonctions d'une courbe lisse et projective sur $k$,
toujours avec $X$ courbe elliptique. Ils
ont formul\'e la conjecture analytique et l'ont reli\'ee avec la
formulation Taub\'erienne,
montrant le Th\'eor\`eme
\ref{thmA} dans ce cas. Plus tard, Wazir dans \cite[Theorem 1.1]{wa}
a g\'en\'eralis\'e le r\'esultat
de Rosen-Silverman au cas d'une vari\'et\'e $\xx$ de dimension 3
fibr\'ee au dessus d'une surface en
courbes elliptiques.
On peut s'attendre \`a ce que l'\'enonc\'e du th\'eor\`eme reste vrai
pour toute fibration en courbe de genre au
moins un $f:\xx\rightarrow T$ avec $n=\dim\xx=\dim T+1$ mais nous
nous contenterons, dans ce texte,
de traiter le cas des
surfaces.

\end{remark}

\section{Conjecture de Birch et Swinnerton Dyer et conjecture de Tate}

Ce paragraphe contient les motivations pour le th\'eor\`eme principal
ainsi que des consid\'erations de
nature sp\'eculative; on explique en particulier pourquoi il est
n\'ecessaire de remplacer la trace
moyenne
$A_{\mathfrak{p}}(\mathcal{X})$ par la trace moyenne ``r\'eduite"
$A_{\mathfrak{p}}^*(\mathcal{X})$ pour esp\'erer
un
\'enonc\'e du type (\ref{NagaoAr}).

\begin{notation}Comme il est d'usage, pour deux fonctions $f(t)$ et
$g(t)$  d\'efinies au voisinage de $\omega$, on \'ecrira,
$f(t)=O(g(t))$, si il existe un voisinage de $\omega$ et
une constante $C$ tels que sur ce voisinage $|f(t)|\le Cg(t)$. On
\'ecrira $f(t)=o(g(t))$ si
$\lim_{t\to\omega}f(t)/g(t)=0$.
\end{notation}

Si $\mathcal{X}$ est une vari\'et\'e lisse projective, d\'efinie sur
un corps de nombres $k$, ayant bonne
r\'eduction hors d'un ensemble fini  de places $S$, on notera, en
n\'egligeant un nombre fini de
facteurs Euleriens~:
$$L_2(\mathcal{X}/k,s)=\prod_{\mathfrak{p}\notin S}
\det(1-F_{\mathfrak{p}}q_{\mathfrak{p}}^{-s}\,|\,
H^2_{\text{\'et}}(\mathcal{X}\times_k\overline{k},\mathbb{Q}_{\ell}))^
{-1}.$$ Pour $\mathfrak{p}$ ne divisant pas
le conducteur, les valeurs propres du Frobenius $F_{\pp}$ sont de
modules $\qp$. Le produit d'Euler
converge pour
$\Re(s)>2$. La Conjecture \ref{conjT} donne une interpr\'etation du
rang du groupe de N\'eron-Severi
$\text{NS}(\mathcal{X})$ en termes de cette fonction. La Conjecture
\ref{conjT} est en fait formul\'ee pour
une vari\'et\'e de dimension quelconque \cite[5.5]{ram}.

Si $A$ est une vari\'et\'e ab\'elienne d\'efinie sur un corps de
nombres $k$, on notera simplement
$$L(A/k,s)=\prod_{\mathfrak{p}}
\det(1-F_{\mathfrak{p}}q_{\mathfrak{p}}^{-s}\,|\,
H^1_{\text{\'et}}(A\times_k\overline{k},\mathbb{Q}_{\ell})^{I_{\mathfrak{p}} })^{-1},$$ sa fonction $L$. Pour
$\mathfrak{p}$ ne divisant pas le conducteur, les valeurs propres du Frobenius $F_{\pp}$, disons
$\alpha_{\mathfrak{p},j}$, sont de modules $q_{\pp}^{1/2}$. Le produit d'Euler converge pour $\Re(s)>3/2$. La
conjecture centrale de la th\'eorie donne une interpr\'etation du rang du groupe de Mordell-Weil $A(k)$ en
termes de
cette fonction.

\begin{conjecture}[Conjecture de Birch et Swinnerton-Dyer
\cite{ta2}]\label{conjBSD} Soit $A$ une vari\'et\'e ab\'elienne
d\'efinie sur un corps de nombres $k$. La
fonction
$L(A/k,s)$ admet un prolongement analytique \`a $\mathbb{C}$ et
poss\`ede en $s=1$ un z\'ero d'ordre
$\rg(A(k))$.
\end{conjecture}


\begin{remark}Il s'agit ici de la premi\`ere partie de la conjecture,
nous ne discuterons pas  la seconde partie qui d\'ecrit le
coefficient dominant de $L(A/k,s)$ en
$s=1$.
\end{remark}


Un calcul simple et classique montre que
\begin{align*} -\frac d{ds}(\log(L(A/k,s)))&=\sum_{\mathfrak{p},m\geq 1}\sum_j
\alpha_{\mathfrak{p},j}^m\log(\qp)\qp^{-ms}\\
&=\sum_{\mathfrak{p}}a_{\mathfrak{p}}(A)\log(\qp)\qp^{-s}+h_1(s)\\
&=\sum_{\mathfrak{p}}a_{\mathfrak{p}}(A)\log(\qp)\qp^{-s}+
\sum_{\mathfrak{p}}\sum_j\alpha_{\mathfrak{p},j}^2\log(\qp)\qp^{-2s}+h_2(s)\ \
\end{align*} avec $h_1(s)$ holomorphe sur $\Re(s)>1$ et $h_2(s)$
holomorphe sur $\Re(s)>5/6$. D\'efinis\-sons
$M(A/k,s)=\sum_{\mathfrak{p}}-a_{\mathfrak{p}}(A)\log(\qp)\qp^{-s}$
ainsi que  la fonction arithm\'etique correspondante $M_A(T):=\sum_{\qp\leq
T}-a_{\mathfrak{p}}(A)\log(\qp)$; alors le
prolongement analytique de $L(A/k,s)$ jusqu'\`a $\Re(s)=1$ \'equivaut
\`a celui de $M(A/k,s)$ et un z\'ero d'ordre
$r$ en $s=1$ pour $L(A/k,s)$ \'equivaut \`a un p\^ole simple pour
$M(A/k,s)$ en $s=1$ avec r\'esidu
$r$. Remarquons aussi
que l'holomorphie de $M(A/k,s)$ pour $\Re(s)>1$ \'equivaut \`a
l'holomorphie de $L'(A/k,s)/L(A/k,s)$ pour $\Re(s)>1$
et donc \`a l'hypoth\`ese de Riemann g\'en\'eralis\'ee pour
$L(A/k,s)$ ou encore, par des arguments
analytiques classiques
\`a la propri\'et\'e $M_A(T)=O_{\epsilon}\left(T^{1+\epsilon}\right)$
pour tout $\epsilon>0$. De plus il
semble raisonnable de conjecturer que $h_1(s)$ devrait \^etre
holomorphe sur $\Re(s)=1$; cela
proviendrait de l'holomorphie et
la non annulation sur $\Re(s)=2$ de la fonction
$\prod_{\mathfrak{p},j}(1-\alpha_{\mathfrak{p},j}^2\qp^{-s})^{-1}$,
ce qui est d\'emontr\'e dans   le cas o\`u $A$ est modulaire (par la
m\'ethode de Rankin,
compl\'et\'ee par
   Shimura).

Acceptant cela, on peut alors s'interroger sur la vraisemblance de
l'\'equivalence $M_A(T)\sim rT$; en effet,
d'une part, si cette \'equivalence \'etait vraie alors on en
d\'eduirait facilement que $M(A/k,s)$
aurait un  p\^ole simple
en $s=1$ avec r\'esidu $r$ et donc la conjecture de Birch et
Swinnerton-Dyer serait v\'erifi\'ee pour $A/k$,
si
$r=\rg(A(k))$, d'autre part, m\^eme en supposant l'hypoth\`ese de
Riemann g\'en\'eralis\'ee et la conjecture de
Birch et  Swinnerton-Dyer pour $L(A/k,s)$, on ne peut en d\'eduire
par les arguments Taub\'eriens
classiques l'\'equivalence
$M_A(T)\sim rT$ \`a cause de la pr\'esence d'une infinit\'e de
p\^oles pour $M(A/k,s)$ (resp. de z\'eros
pour
$L(A/k,s)$) sur la droite $\Re(s)=1$. La pr\'esence d'une infinit\'e
de p\^oles rend m\^eme fort peu
vraisemblable l'\'equivalence $M_A(T)\sim rT$. On conclut ces
consid\'erations en notant qu'il
semble donc illusoire d'esp\'erer
remplacer la conclusion du th\'eor\`eme principal par
$$\lim_{T\to\infty}\frac 1T\sum_{\qp\leq
T}-A_{\mathfrak{p}}(\mathcal{X})\log\qp={\rm rang}\,\left(J_X
(K)\right)??$$ (sauf, bien s\^ur, dans le cas  o\`u
$B$ est nulle). Par contre la version analytique
$$
\res_{s=1}\left(\sum_{\mathfrak{p}\notin
S}-A_{\mathfrak{p}}(\mathcal{X})\frac{\log(q_{\mathfrak{p}})}
{q_{\mathfrak{p}}^s}\right)=\text{rang}\left(J_X(K)\right)?
$$ d\'ecoulerait des conjectures de Tate pour $\mathcal{X}$ et Birch
et Swinnerton-Dyer pour $B$. En effet, si l'on admet le prolongement
analytique, alors le r\'esidu en $s=1$ est
\'egal
\`a:
\begin{multline*}\rg\left(J_X(K)\right)+\left(\ord_{s=1}(L(B/k,s))-\text{rang}(B(k))\right)\\
-\left(\rg\left(\ns(\xx/k)\right)+\ord_{s=2}(L_2(\xx/ k,
s))\right).\end{multline*}

Ainsi la situation est plus favorable sur les corps de fonctions;
cela est d\^u au fait que, du point de vue de
la th\'eorie analytique des nombres, on effectue une double moyenne
et que l'on doit consid\'erer le
comportement d'une
s\'erie de Dirichlet au bord de son demi-plan de convergence, et non
\`a l'int\'erieur de la bande critique.

\begin{remark} Il serait peut-\^etre plus naturel, au lieu de
d\'efinir la moyenne des traces de Frobenius $A_{\pp}(\mathcal{X})$
comme nous l'avons fait (imitant en
cela Nagao~\cite{nag2} et Rosen-Silverman~\cite{rosi}), de poser

$$A'_{\pp}(\mathcal{X})=
\frac 1{\#C_{\pp}(\kappap)}\sum_{c\in
C_{\pp}(\kappap)}a_{\pp}(\mathcal{X}_{c}) \quad\hbox{et}\quad
A_{\pp}^{'*}(\mathcal{X}):=A'_{\pp}(\mathcal{X})-a_{\pp}(B).$$  On
peut observer qu'une telle modification
ne changerait pas les \'enonc\'es en vue. En effet, on voit
ais\'ement que l'\'egalit\'e voulue
$$\lim_{T\to\infty}\frac 1T\sum_{\substack{\pp\notin S\\ q_{\pp}\le
T}}(A_{\pp}^*(\mathcal{X})-A^{'*}_{\pp}(\mathcal{X}))\log(\qp)=0$$
est v\'erifi\'ee si et seulement
si
$$N(T):=\sum_{\substack{\pp\notin S\\ q_{\pp}\le
T}}a_{\pp}(C)A^*_{\pp}(\mathcal{X})\frac{\log(\qp)}{ \qp}=o(T).
$$ Cette derni\`ere estimation est donc vraie si $C=\ppb^1$.
Observons maintenant que $A_{\pp}^*(\xx)$, qui est trivialement
$O(q_{\pp}^{1/2})$ (comme moyenne des
$a_{\pp}(\xx_c)$ qui sont $O(q_{\pp}^{1/2})$) est en fait un $O(1)$
\`a cause de la formule
(\ref{form5}) (d\'emontr\'ee plus loin au
paragraphe 4, o\`u chaque terme est d\'efini) qui indique que
$$A_{\pp}^*(\xx)=-\frac{b_{\pp}(\xx)}{\qp}+\tr(F_{\pp}|\mathcal{S})+O(
\qp^{-1/2})$$ Or $b_p(\xx)$ est la trace d'un op\'erateur avec
valeurs propres de module $\qp$ donc
$b_{\pp}(\xx)=O(\qp)$, alors que $F_{\pp}$ agit sur $\mathcal{S}$
avec pour valeurs propres des racines de
l'unit\'e (i.e., de fa\c con quasi-unipotente), donc on a bien
$A_{\pp}^*(\xx)=O(1)$. Ainsi on
obtient
$$N(T)=\sum_{\qp\leq T}a_{\pp}(C)A_{\pp}^*(\xx)\frac{\log
\qp}{\qp}\ll\sum_{\qp\leq T}|a_{\pp}(C)|\frac{\log \qp}{\qp}
\ll\sum_{\qp\leq T}\frac{\log \qp}{\sqrt{\qp}}\ll
\log T\,\sqrt{T}.$$
\end{remark}

\section{Outils g\'eom\'etriques}

Tout comme dans Rosen-Silverman \cite{rosi}, l'id\'ee de la preuve
est de compter le nombre des points rationnels
de
$\mathcal{X}_{\mathfrak{p}}$ sur $\kappa_{\mathfrak{p}}$ de deux
mani\`eres. La premi\`ere consiste \`a employer
la formule de Lefschetz pour $\mathcal{X}_{\mathfrak{p}}$, la
deuxi\`eme \`a compter le nombre des
points dans chaque
fibre $\mathcal{X}_{\mathfrak{p},c}$ pour $c\in
C_{\mathfrak{p}}(\kappa_{\mathfrak{p}})$ et calculer la somme.
Pour obtenir ces deux r\'esultats on a besoin d'une part d'un lemme
g\'eom\'etrique (Lemme
\ref{geom1}) qui permet de
d\'eterminer la cardinalit\'e de
$\mathcal{X}_{\mathfrak{p},c}(\kappa_{\mathfrak{p}})$ et d'autre part
d'une
formule (Proposition \ref{geom3}) analogue \`a celle de Shioda-Tate
pour les surfaces elliptiques,
g\'en\'eralis\'ee par
Shioda aux fibrations de genre sup\'erieur admettant une section
\cite[Theorem 1]{sh1}. On d\'eveloppe ces
r\'esultats dans ce paragraphe. Pour obtenir le premier, nous
g\'en\'eralisons les calculs
explicites au cas par cas de
\cite{rosi,wa} par le Lemme \ref{geom1}. Enfin nous incluons une
preuve de la formule de Shioda-Tate  \cite{sh1}
car, d'une part, nous ne souhaitons pas faire l'hypoth\`ese que la
fibration poss\`ede une section,
d'autre part, le fait de
travailler en caract\'eristique z\'ero simplifie la preuve (les
sch\'emas en groupes, notamment $\pic^0$,
sont r\'eduits).

Nous commen\c{c}ons par rappeler le lemme suivant qui indique
quelques propri\'et\'es des fibres singuli\`eres
de familles de courbes (cf. \cite{arwi}).

\begin{lemma}\label{lem1}\cite[(1.1)]{arwi} Soit $R$ un anneau de
valuation discr\`ete avec corps r\'esiduel alg\'ebriquement clos
$\kappa$. Soit $\mathcal{V}$ un sch\'ema r\'egulier
de dimension 2 propre sur \linebreak $\spec(R)$. On suppose que la
fibre ferm\'ee ${V}$ soit connexe
et de dimension 1. En
tant que sous-sch\'ema de $\mathcal{V}$, $V$ est un diviseur de
Cartier effectif que l'on \'ecrit sous la
forme
$V=\sum_{i=1}^nr_iC_i$. La matrice d'intersection $||(C_i.C_j)||$
d\'efinit une forme quadratique dont le noyau
est engendr\'e par $V$ et qui est d\'efinie n\'egative sur
$(\bigoplus_{i=1}^n\qq C_i)/\qq V$.
\end{lemma}

Un point clef dans notre argument est le lemme suivant. Nous
remercions le {\it referee} pour nous avoir
indiqu\'e la preuve de ce lemme (nous proposions une preuve beaucoup
plus compliqu\'ee, bas\'e sur
l'analyse des fibres sp\'eciales singuli\`eres donn\'ee par Artin et
Winters \cite{arwi}, qui de plus
ne donnait le r\'esultat qu'\`a $O(q_{\pp}^{1/2})$ pr\`es).

\begin{lemma}\label{geom1}Soient $\mathfrak{p}\notin S$, $c\in
C_{\mathfrak{p}}(\kappa_{\mathfrak{p}})$,
$\mathcal{X}_{\mathfrak{p},c}=f_{\mathfrak{p}}^{-1}(c)$ et $m_c$ le
nombre des composantes $\kappa_{\mathfrak{p}}$-rationnelles de
$\mathcal{X}_{\mathfrak{p},c}$. Le
nombre de points
$\kappa_{\pp}$-rationnels de $\xx_{\pp,c}$ s'exprime comme suivant:
$$
\#\xx_{\pp,c}(\kappa_{\pp})=1-a_{\pp}(\xx_{\pp,c})+q_{\pp}m_c.
$$

\end{lemma}

\begin{proof}
Par le th\'eor\`eme de Weil, pour tout
$c\in(C_{\pp}-\Delta_{\pp})(\kappa_{\pp})$ on a
$\#\xx_{\pp,c}(\kappa_{\pp})=q_{\pp}+1-a_{\pp}(\xx_{\pp,c})$ et $m_c=1$.

Supposons $c\in\Delta_{\pp}(\kappa_{\pp})$. On peut passer en premier lieu \`a la structure r\'eduite de $\xx_{\pp,c}$,
en effet cela ne change  ni le d\'ecompte des points rationnels ni la cohomologie. Par la Formule de Lefschetz (cf.
\cite[Rapport, Corollaire 5.4]{sga4demi} et \cite[(1.5)]{del74}, pour tout $c\in C_{\pp}(\kappa_{\pp})$ on a
$\#\xx_{\pp,c}(\kappa_{\pp})=\sum_{i=0}^2(-1)^i\tr(\ov{F}_{\pp}\,|\,H^ i_{\unc}(\ov{\xx}_{\pp,c},\qq_{\ell}))$. Ici
$H^i_{\unc}$ d\'esigne la cohomologie $\ell$-adique \`a support propre. Comme $\xx_{\pp,c}$ est connexe,
$\tr(\ov{F}_{\pp}\,|\,H^0_{\unc}(\ov{\xx}_{\pp,c},\qq_{\ell}))=1$. De plus, par d\'efinition,
$a_{\pp}(\xx_{\pp,c})=\tr(\ov{F}_{\pp}\,|\,H^1_{\unc}(\ov{\xx}_{\pp,c}$ \linebreak $,\qq_{\ell}))$.

Soit $\vartheta_{\pp}:\cac_{\pp}\to\ov{\xx}_{\pp,c}$ une normalisation de $\ov{\xx}_{\pp,c}$. La courbe $\cac_{\pp}$
est une r\'eunion disjointe de courbes projectives lisses. Le morphisme ${\vartheta}_{\pp}$ induit une suite
spectrale
de Leray $E^{r,s}:=H^r_{\unc}(\cac_{\pp},R^s({\vartheta}_{\pp})_*\qq_{\ell})\Rightarrow
E^{r+s}:=H^{r+s}_{\unc}(\ov{\xx}_{\pp,c},\qq_{\ell})$. Elle induit la suite exacte suivante:

\begin{align*}\begin{matrix}
0&\to&H^1_{\unc}(\cac_{\pp},({\vartheta}_{\pp})_*\qq_{\ell})&\to&H^1_{
\unc}(\ov{\xx}_{\pp,c},\qq_{\ell})&\to&
H^0_{\unc}(\cac_{\pp},R^1_c(({\vartheta}_{\pp})_*\qq_{\ell}))\\
&\to&H^2_{\unc}(\cac_{\pp},({\vartheta}_{\pp})_*\qq_{\ell})&\to&H^2_{\
unc}(\ov{\xx}_{\pp,c},\qq_{\ell})&\to&
H^1_{\unc}(\cac_{\pp},R^1_{\unc}(({\vartheta}_{\pp})_*\qq_{\ell})).
\end{matrix}\end{align*}

Il suit de \cite[Theorem 3.2 (b), Chapter VI, p. 228]{mi2} que pour tout $i>0$ on a
$R^i_{\unc}({\vartheta}_{\pp})_*=0$, donc $H^2_{\unc}(\cac_{\pp},({\vartheta}_{\pp})_*\qq_{\ell})\cong
H^2_{\unc}(\ov{\xx}_{\pp,c},\qq_{\ell})$. De plus, le premier groupe est une somme directe de $\qq_{\ell}(-1)$. Pour
chaque composante $\kappa_{\pp}$-rationnelle de $\cac_{\pp}$, l'automorphisme de Frobenius $\ov{F}_{\pp}$ agit sur la
composante correspondante dans $H^2_{\unc}(\cac_{\pp},({\vartheta}_{\pp})_*\qq_{\ell})$ par multiplication par
$q_{\pp}$. Les autres composantes de $\cac_{\pp}$ sont permut\'ees entre elles, a fortiori la trace de $\ov{F}_{\pp}$
dans les composantes correspondantes dans $H^2_{\unc}(\cac_{\pp},({\vartheta}_{\pp})_*\qq_{\ell})$ est \'egale \`a 0.
D'o\`u $\tr(\ov{F}_{\pp}\,|\,H^2_{\unc}(\ov{\xx}_{\pp,c},\qq_{\ell}))=q_{\pp} m_c$ et cela ach\`eve la preuve du lemme.

\end{proof}

Nous passons maintenant \`a  l'\'etude de la $K/k$-trace not\'ee $B$.

Soit $\alb(\xx)$, respectivement $\alb(C)$, la vari\'et\'e Albanese
de $\xx$, respectivement de $C$. Soit $(\tau,B)$
la
$K/k$-trace de $J_X$. Soit $\pic^0(C)$ le groupe des diviseurs
alg\'ebrique\-ment \'equivalents \`a z\'ero de $C$
et
$f^*:\pic^0(C)\to\text{Pic}^0(\mathcal{X})$ l'application
``pull-back'' d\'eduite de $f$.

   Soit
$\imath:X\to\mathcal{X}$ l'inclusion de la fibre g\'en\'erique et
$\imath^*:{\rm Div}(\mathcal{X})\to\text{Div}(X)$
la restriction des diviseurs \`a la fibre g\'en\'erique $X$. Par
construction de la $K/k$-trace,
l'homomorphisme induit de
$\pic^0(\xx)$ vers $\pic^0(X)$ se factorise par un homomorphisme
$b:\pic^0(\xx)\to B$.

\begin{proposition}\label{geom4}
Le groupe $\kker(f^*)$ est fini et
est m\^eme trivial si $f$ admet une section.  La suite  de
vari\'et\'es ab\'eliennes
\begin{equation}\label{2.9}
0\to\pic^0(C)\overset{f^*}\longrightarrow\pic^0(\xx)\overset{b}\longrightarrow B\to 0
\end{equation}
est exacte si $f$ admet une
section et exacte \`a des groupes finis pr\`es en g\'en\'eral. En particulier,
$H^1_{\et}(\xbar,\qq_{\ell})\cong H^1_{\et}(\overline{C},\qq_{\ell})\oplus
H^1_{\et}(\overline{B},\qq_{\ell})$
   A fortiori, on a
pour presque tout id\'eal premier $\pp$ de $\mathcal{O}_k$,
$a_{\pp}(\xx)=a_{\pp}(C)+a_{\pp}(B)$.
\end{proposition}

\begin{remark}
La formule $B=\pic^0(\xx)/f^*\pic^0(C)$ est donn\'ee
(sans preuve) dans \cite{ta1}
   et
l'exactitude de la suite  est d\'emontr\'ee par Shioda dans
\cite{sh1} sous l'hy\-po\-th\`e\-se d'existence
d'une section mais en compl\`ete g\'en\'eralit\'e concernant la
caract\'eristique et la structure
de sch\'ema (le r\'esultat
est attribu\'e \`a Raynaud).
\end{remark}

\begin{proof}
Il est clair que $\text{Im}(f^*)\subset\kker(b)$, car
$b\circ f^*=0$. Si $D$ est un diviseur sur $\xx$ dont la classe est
dans $\kker(b)$, on voit que
$D=\text{div}(g)+V$, o\`u $\text{div}(g)$ est le diviseur d'une
fonction rationnelle
$g\in\overline{k}(\xx)$ et $V$ est un diviseur
vertical, i.e., \`a support dans les fibres de $f$. Mais $D$ et donc
$V$ sont alg\'ebriquement \'equivalents
\`a z\'ero, et donc num\'eriquement triviaux, donc, d'apr\`es le
Lemme \ref{lem1}, $V$ est une
somme de fibres,
c'est-\`a-dire qu'il est de la forme $f^*(D')$. {\it A priori} on
doit prendre $D'$ dans $\pic^0(C)\otimes\qq$ et
donc il existe $m\geq 1$ tel que $mV\in f^*\pic^0(C)$, mais si $f$
admet une section, il n'y a pas
de fibres multiples et on
peut prendre $m=1$; en g\'en\'eral $m$ est born\'e par la
multiplicit\'e des fibres. On a donc bien
$\text{Im}(f^*)=\kker(b)$.

Soit $C_0\subset \mathcal{X}$ une courbe irr\'eductible ferm\'ee
telle que $f_{|C_0}:C_0\to C$ soit dominante, donc
un morphisme fini de degr\'e disons $d$. Notons $j$ l'inclusion
$C_0\hookrightarrow\xx$, donc
$f_{|C_0}=f'=f\circ j$.
L'hypoth\`ese $f^*(z)=0$ entra\^{\i}ne ${f'}^*(z)=0$ donc
$dz=f'_*{f'}^*(z)=0$. Ainsi $\kker(f^*)\subset
\pic^0(C)[d]$ est fini (et m\^eme trivial si on peut choisir $d=1$,
i.e., si $C_0$ est une section
de $f$).

Pour prouver que $b$ est surjective, il suffit de prouver que $\dim(B)\leq \dim(\pic^0(\mathcal{X}))$ \linebreak
$-\dim(\pic^0(C))$, ou encore que $\dim(\alb(C))+\dim(B)\leq\dim(\alb(\mathcal{X}))$ (il y aura alors \'egalit\'e).
L'homomorphisme $\tau:B\to J_X$ est injectif, parce que $k$ est un corps de nombres. Son homomorphisme dual Soit
$\tau^{\vee}:J_X\to B^{\vee}$ c'est en fait la $K/k$-image de $J_X$, cf. \cite[Chapter VIII]{la2}. Notons que par
extension  des scalaires, et en tenant compte que $\ov{K}\supset\ov{k}(C)$, la $K/k$-trace et la $K/k$-image de $J_X$
co\"{\i}cident avec la $\ov{K}/\ov{k}$-trace et la $\ov{K}/\ov{k}$-image de $J_X$. Choisissons un morphisme fini
$X\rightarrow J_X$ d\'efini sur $K$ et dont l'image engendre $J_X$. En composant avec $\tau$, on obtient une
application qui se prolonge en une application rationnelle $\mathcal{X}\dashrightarrow B^{\vee}\times C$ d\'efinie sur
$k$. En composant avec l'application canonique $\pi_C:C\rightarrow\alb(C)$, on obtient une application, \textit{a
priori} rationnelle $\alpha:\mathcal{X}\dashrightarrow B^{\vee}\times\alb(C)$ qui doit en fait \^etre un morphisme,
puisque $\xx$ est lisse (cf. \cite[Chapter II, Theorem 1]{la2}). Ce dernier morphisme induit donc, par la propri\'et\'e
universelle de la vari\'et\'e d'Albanese, un homomorphisme $\hat{\alpha}:\alb(\mathcal{X})\rightarrow
B^{\vee}\times\alb(C)$ dont on montre qu'il est surjectif,  ce qui ach\`eve la d\'emonstration de l'exactitude. Pour
v\'erifier que $\hat{\alpha}$ est surjective (donc c'est une isog\'enie), il faut montrer que $\alpha(\mathcal{X})$
engendre la vari\'et\'e ab\'elienne $B^{\vee}\times\alb(C)$. L'image de la fibre g\'en\'erique $X$ dans sa Jacobienne
$J_X$ engendre $J_X$, donc son image par $\tau^{\vee}$ engendre $B^{\vee}$ et comme l'image $\pi_C(C)$ engendre
$\alb(C)$, on peut conclure.

Pour l'isomorphisme entre les groupes de cohomologie \'etale, comme $\kker(f^*)$ est fini, la tensorisation par
$\qq_{\ell}$ donne $H^1_{\text{\'et}}(\pic^0(C),\qq_{\ell})\cong H^1_{\text{\'et}}(\overline{C},\qq_{\ell})$, donc
l'affirma\-tion sur la trace est imm\'ediate.
\end{proof}

Nous passons maintenant \`a la formule de Shioda-Tate.

L'ensemble $I_X=\{(D.X)\,|\,D\in{\rm Div}(\mathcal{X})\}$ est un
id\'eal de $\mathbb{Z}$. Soit $(D_0.X)$
un g\'en\'erateur de $I_X$, o\`u $D_0$ est un diviseur sur $\mathcal{X}$. Soit
$\imath:X\to\mathcal{X}$ l'inclusion de la
fibre g\'en\'erique et $\imath^*:{\rm
Div}(\mathcal{X})\to\text{Div}(X)$ la restriction des diviseurs \`a la
fibre g\'en\'erique $X$. D\'efinissons
$\psi:\text{Pic}(\mathcal{X})\to J_X$ par
$$\psi(\text{cl}(D))=\imath^*\left(D-\frac{(D.{X})}{(D_0.{X})}D_0\right).$$ On obtient alors
d'apr\`es le Lemme ci-dessous (\ref{geom2}) une suite exacte
\begin{equation}\label{2.3}
0\to\kker(\psi)\to\pic(\xx)\overset{\psi}\longrightarrow J_X\to0.
\end{equation}

\begin{definition}Soit $\tilde{\ssc}$ (respectivement $\ssc$) le sous-groupe de
$\pic(\xx)$ (respectivement de $\ns(\xx)$) engendr\'e par les classes
de $D_0$  et des composantes  des fibres  de
$f$. On notera respectivement $\tilde{\ssc}_{\qq}$ et $\ssc_{\qq}$
les espaces vectoriels obtenus
par tensorisation par
$\qq$.
\end{definition}

\begin{remark}
Soit $\pp\notin S$  et $X_{\pp}$ la fibre g\'en\'erique de $f_{\pp}$.
L'ensemble $I_{X_{\pp}}:=\{(D\cdot
X_{\pp})\;|\; D\in\ddiv(\xx_{\pp})\}$ est un id\'eal de $\zz$. Soit
$\mathcal{D}_{1,\pp}$ un
diviseur de $\xx_{\pp}$ tel que
$(\mathcal{D}_{1,\pp}\cdot X_{\pp})$ engendre $I_{X_{\pp}}$. Soit
$\ssc_{\pp}$ le sous-espace de
$\ns(\xx_{\pp})$ engendr\'e par les classes de $\mathcal{D}_{1,\pp}$
et des composantes des fibres
de $f_{\pp}$. Pour tout $\pp$ (sauf
un nombre fini), $\ssc_{\pp}$ co\"{\i}ncide avec la r\'eduction de
$\ssc$ modulo $\pp$. Donc, apr\`es avoir \'elargi
si n\'ecessaire $S$, on peut supposer cela vrai pour tout $\pp$ hors
de $S$. En particulier, pour
$\pp$ hors de $S$, on
aura $\tr(F_{\pp}|\ssc^{I_{\pp}})=\tr(\overline{F}_{\pp}|\ssc_{\pp})$.
\end{remark}

\begin{lemma}\label{geom2}L'application $\psi$ est surjective et son
noyau $\kker(\psi)$ est \'egal \`a $\tilde{\ssc}$. Une base de
${\ssc}_{\qq}$ est fournie par
   les classes de $D_0$, d'une fibre lisse $\mathcal{F}$ de $f$ et par
les composantes, sauf une,
    de chaque fibre singuli\`ere non irr\'eductible
$f^{-1}(c)=\mathcal{X}_c$ de $f$. Soit $n_c$ le nombre des
composantes de $\mathcal{X}_c$. On a en
particulier,
$$\rg(\ssc)=2+\sum_{c\in C}(n_c-1).$$
\end{lemma}

\begin{proof} Soit $P\in J_X$ tel que $P$ repr\'esente la classe
$\text{cl}(\mathcal{D})$ d'un diviseur
$\mathcal{D}\in\text{Div}({X})$. Soit $\overline{\mathcal{D}}$ la
cl\^oture
de Zariski de $\mathcal{D}$ dans $\mathcal{X}$. C'est un diviseur de
Weil (donc de Cartier, puisque
$\mathcal{X}$ est
lisse) et il est imm\'ediat qu'on a
$\psi(\text{cl}(\overline{\mathcal{D}}))=P$, o\`u
$\text{cl}(\overline{\mathcal{D}})$ d\'enote la classe de
$\overline{\mathcal{D}}$. Ainsi $\psi$ est bien surjective.

Il est clair que $\tilde{\ssc}\subset\kker(\psi)$. Supposons que
$\text{cl}(D)\in\kker(\psi)$. Alors, il existe
une fonction $x\in \overline{K}(X)$ telle que
$\text{div}(x)=\imath^*(D-{((D.X)/(D_0.X))}D_0)$.
Soit $\tilde{x}$ une
fonction rationnelle sur $\mathcal{X}$ telle que $\tilde{x}_{|X}=x$.
En particulier,
$\imath^*\text{div}(\tilde{x})={\text{div}}(x)$. Il est clair que, au
niveau des diviseurs,
   $\kker(\imath^*)$ est
engendr\'e par  les classes des composantes irr\'eductibles
   des fibres
${\mathcal{X}}_c$ de $f$. Donc
$\text{div}(\tilde{x})-(D-{((D.X)/(D_0.X))}D_0)$ est une somme des
composantes de
ces fibres. On en tire bien que $\text{cl}(D)\in\tilde{\ssc}$ et, par suite
$\kker(\psi)=\tilde{\ssc}$.

Pour la derni\`ere affirmation, il est clair que les classes de
l'\'enonc\'e forment un syst\`eme g\'en\'erateur
car toutes les fibres sont alg\'ebriquement \'equivalentes. Le fait
qu'elles soient ind\'ependantes
se v\'erifie en
utilisant la th\'eorie d'intersection et le Lemme \ref{lem1}. La
formule donnant le rang de $\ssc$ est
alors imm\'ediate.
\end{proof}

\begin{proposition}[Formule de Shioda-Tate]\label{geom3}  Il existe
un isomorphisme de
\newline $\qq[\gk]$-modules
$$\ns(\xx)\otimes\qq\cong\left(\left(\frac{J_X(\ov{k}(C))}{\tau
B(\ov{k})}\right)\otimes\qq\right)\oplus\ssc_{\qq}.$$ En particulier, on a
\begin{equation}\label{2.11}\rg(\ns(\xx/k))=\rg\left(\frac{J_X(K)}{\tau
B(k)}\right)+\rg(\ssc^{\gk}).
\end{equation}

\end{proposition}

\begin{proof} Remarquons que, si l'on travaille sur les $\qq$-espaces
vectoriels, on peut supposer que $D_0$ soit d\'efini sur $k$ (la
seule diff\'erence est que $(D\cdot X)/(D_0\cdot
X)$ sera maintenant un nombre rationnel). Les suites exactes
(\ref{2.9}) et  (\ref{2.3}) et le fait
que
$\kker(\psi)=\tilde{\ssc}$ nous donnent une suite exacte et scind\'ee
apr\`es tensorisation par
$\qq$
\begin{equation}\label{2.10}
0\to\ssc\to\ns(\xx)\overset{\vartheta}\longrightarrow\frac{J_X(\ov{k}(C))}{\tau
B(\ov{k})}\to0,
\end{equation}
le r\'esultat en d\'ecoule.
\end{proof}

\section{Preuve du Th\'eor\`eme \ref{thmA}}

Soient $\mathcal{S}_{\ell}(1)=\mathcal{S}\otimes
T_{\ell}(\overline{k})$ le ``twist'' \`a la Tate  du
$\gk$-module
$\mathcal{S}$ et $L(\mathcal{S}_{\ell}(1),s)$ la fonction $L$ d'Artin
associ\'ee \`a la repr\' esentation donn\'ee
par
$\mathcal{S}_{\ell}(1)$. Comme dans \cite{rosi}, la Proposition
\ref{geom3} nous sugg\`ere de consid\'erer la
fonction
$N(\mathcal{X}/k,s)={L_2(\mathcal{X}/k,s)}/$ \linebreak
${L(\mathcal{S}_{\ell}(1),s)}$.

\begin{proposition}\label{anal1a} Avec les notations pr\'ec\'edentes, on a
\begin{equation}\label{anal1}\rg\left(\frac{J_X(K)}{\tau B(k)}\right)
+\ord_{s=2}(N(\mathcal{X}/k,s))=\rg(\ns(\xx/k))+\ord_{s=2}
(L_2(\mathcal{X}/k,s)).\end{equation}
\end{proposition}

\begin{proof} D'apr\`es la Proposition \ref{geom3} et la d\'efinition
de $N(\mathcal{X}/k,s)$, il suffit de prouver que
$$-\ord_{s=2}(L(\mathcal{S}_{\ell}(1),s))=\rg(\mathcal{S}^{\gk}).$$
Cela est une cons\'equence du r\'esultat suivant.
\end{proof}

\begin{proposition}[Artin, Brauer]\label{arbr}
\cite[Proposition 1.5.1]{rosi} Soit $V$ un $\mathbb{Q}$-espace
vecto\-riel de dimension finie avec une action
continue de $\gk$, et soit $L(V,s)$ la fonction $L$ d'Artin associ\'ee.
\begin{enumerate}
\item $L(V,s)$ a une continuation m\'eromorphe \`a $\mathbb{C}$.
\item $L(V,s)$ est holomorphe sur la droite
$\Re(s)=1$
\`a l'exception \'eventuelle du point $s=1$, o\`u l'on a
$\ord_{s=1}(L(V,s))=-\dim(V^{\gk})$. \item $L(V,s)$
ne s'annule pas sur la droite $\Re(s)=1$.
\end{enumerate}
\end{proposition}

\begin{proposition}\label{anal2a} On a la formule
\begin{equation}\label{anal2}\frac
d{ds}(\log(N(\mathcal{X}/k,s)))=\sum_{\pp\notin
S}A_{\mathfrak{p}}^*(\mathcal{X})\frac{\log(q_{\mathfrak{p}})}
{q_{\mathfrak{p}}^{s-1}}+h(s),\end{equation} avec $h(s)$ fonction
holomorphe sur le demi-plan
$\Re(s)>3/2$.
\end{proposition}

\begin{proof} La fonction-z\^eta de $\mathcal{X}$ s'\'ecrit comme
$$\zeta(\mathcal{X}/k,s)=\prod_{i=0}^4P_i(s)^{(-1)^{i+1}},$$ o\`u
$P_i(s)=\prod_{\pp}\det(1-F_{\mathfrak{p}}q_{\mathfrak{p}}^{-s}\,|\,H^
i_{\text{\'et}}(\overline{\mathcal{X}},
\mathbb{Q}_{\ell})^{I_{\pp}})^{-1}$. Notons que $P_0(s)=\zeta_k(s)$
et $P_4(s)=\zeta_k(s+2)$. En plus
$P_3(s)=P_1(s+1)$
\cite[5.6]{ram}. La formule de Lefschetz et la Proposition
\ref{geom4} impliquent

\begin{equation}\label{form2}
\#\xx_{\pp}(\kappap)=1-a_{\pp}(C)-a_{\pp}(B)+b_{\pp}(\xx)-\qp
a_{\pp}(C)-\qp a_{\pp}(B)+q_{\pp}^2,\end{equation}
o\`u
$b_{\pp}(\xx)$ d\'enote la trace de $F_{\pp}$ dans
$H^2_{\text{\'et}}(\overline{\xx},\qqell)^{I_{\pp}}$. D'autre
part par le Lemme \ref{geom1} on a
\begin{equation}\begin{aligned}\label{form2a}&\#\xx_{\pp}(\kappap)=
\sum_{c\in C_{\pp}(\kappap)}\#\xx_{\pp,c}(\kappap)\\ &=\sum_{c\in
C_{\pp}(\kappap)}(\qp+1+(m_c-1)\qp-a_{\pp}(\xx_{\pp,c}))\\
&=(\qp+1)\#C_{\pp}(\kappap)-\qp A_{\pp}(\xx)+\left(\sum_{c\in
C_{\pp}(\kappap)}(m_c-1)\right)\qp.
\end{aligned}\end{equation} Par ailleurs
$\#C_{\pp}(\kappap)=\qp+1-a_{\pp}(C)$.  En tenant compte aussi que
$a_{\pp}(B)=O({\qp}^{1/2})$, on obtient
de (\ref{form2}) et (\ref{form2a}) que
\begin{equation}\label{form3}\qp
A_{\pp}^*(\xx)=-b_{\pp}(\xx)+2\qp+\left(\sum_{c\in
C_{\pp}(\kappap)}(m_c-1)\right)\qp+O(q_{\pp}^{1/2}).\end{equation} Le Frobenius
$\overline{F}_{\pp}$ fixe les classes de $D_0\pmod{\pp}$  et
$\mathcal{F}\pmod{\pp}$ et permute les
composantes  des fibres singuli\`eres de $f_{\pp}$. Donc, la trace de
$\overline{F}_{\pp}$ sur $\ssc_{\pp}$ est
\'egale
\`a $2$ plus le nombre des composantes $\kappap$-rationnelles, i.e.,
\begin{equation}\label{form4}\tr(\overline{F}_{\pp}\,|\,\mathcal{S}_{\pp})
=\tr(F_{\pp}\,|\,\mathcal{S}^{I_{\pp}})=2+\sum_{c\in
C_{\pp}(\kappap)}(m_c-1).\end{equation} Il suit donc
de (\ref{form3}) et (\ref{form4}) que
\begin{equation}\label{form5}A_{\pp}^*(\xx)=-\frac{b_{\pp}(\xx)}{\qp}+\tr(F_
{\pp}\,|\,\mathcal{S}^{I_{\pp}}) +O(q_{\pp}^{-1/2}).\end{equation}

On calcule maintenant les d\'eriv\'ees logarithmiques de
$L(\ssc_{\ell}(1),s)$ et $L_2(\xx/k,s)$. Pour la premi\`ere
on a
\begin{equation}\label{form6}\begin{aligned}\frac
d{ds}(\log(L(\ssc_{\ell}(1),s)))&=\frac
d{ds}\left(\sum_{\pp}-\log(\det(1-F_{\pp}q_{\pp}^{-s}\,|\,\ssc_{\ell}(
1))^{I _{\pp}})\right)\\
&=\sum_{\pp\notin S}-\tr(F_{\pp}\,|\,\ssc)\frac{\log(\qp)}{q_{\pp}^{s-1}}+f_1(s),
\end{aligned}\end{equation} avec
$f_1(s)$ holomorphe sur le demi-plan $\Re(s)>3/2$. Pour la deuxi\`eme,
\begin{equation}\label{form7}\begin{aligned}\frac
d{ds}(\log(L_2(\xx/k,s)))&=\frac d{ds}\left(\sum_{\pp}-\log(\det(1-F_{\pp}q_{\pp}^{-s}\,|\,H^2_{\text{\'et}}(
\overline{\xx},\qqell))^{I_{\pp}})\right)\\ &=\sum_{\pp\notin
S}-b_{\pp}(\xx)\frac{\log(\qp)}{q_{\pp}^s}+f_2(s),\end{aligned}\end{equation}
avec une fonction $f_2(s)$ holomorphe
pour $\Re(s)>3/2$. En effet, pour ce demi-plan les conjectures de Weil (Deligne) nous donnent
$|b_{\pp}(\xx)|\le\qp\dim(H^2_{\text{\'et}}(\overline{\xx},\qqell))$.

Finalement, par la d\'efinition de $N(\xx/k,s)$, (\ref{form6}) et
(\ref{form7}) on obtient
\begin{align*}\frac d{ds}(\log(N(\xx/k,s)))&=\frac
d{ds}(\log(L_2(\xx/k,s)))-\frac d{ds}(\log(L(\ssc_{\ell}(1),s)))\\
&=\sum_{\pp\notin
S}-b_{\pp}(\xx)\frac{\log(\qp)}{q_{\pp}^s}-\sum_{\pp\notin
S}-\tr(F_{\pp}\,|\,\ssc)\frac{\log(\qp)}{q_{\pp}^{s-1}} +h_1(s)\\
&=\sum_{\pp}A_{\pp}^*(\xx)\frac{\log(\qp)}{q_{\pp}^{s-1}}+h_2(s),
\end{align*} avec $h_2(s)$ holomorphe pour $\Re(s)>3/2$ parce que la
s\'erie $\sum_{\pp\notin
S}\qp^{-1/2}\frac{\log(\qp)}{\qp^{\sigma-1}}$ est convergente pour
$\sigma>3/2$,
o\`u
$\Re(s)=\sigma$.
\end{proof}

\begin{proof}[D\'emonstration du Th\'eor\`eme \ref{thmA}] Les
Propositions \ref{anal1a} et \ref{anal2a} et la Conjecture
\ref{conjT} nous donnent
\begin{align*}&\res_{s=1}\left(\sum_{\pp}-A_{\pp}^*(\xx)\frac{\log(\qp)}{q_{
\pp}^s}\right)=-\res_{s=2}\left(\frac
d{ds}(\log(N(\xx/k,s)))\right)=-\ord_{s=2}(N(\xx/k,s))\\
&=-\rg(\ns(\xx))+\rg\left(\frac{J_X(K)}{\tau B(k)}\right)
-\ord_{s=2}L_2(\xx,s)=\rg\left(\frac{J_X(K)}{\tau B(k)}\right).
\end{align*}

Par la Proposition \ref{arbr}, $L(\ssc_{\ell}(1),s)$ ne s'annule pas
sur la droite $\Re(s)=2$. En supposant de plus
que
$L_2(\xx/k,s)$ ne s'annule pas sur la droite $\Re(s)=2$, on peut
appliquer les th\'eor\`emes
Taub\'eriens standard
\cite[Chapter XV]{la} obtenant
\begin{align*}\label{form8}&\lim_{T\to\infty}\frac
1T\left(\sum_{\substack{\pp\notin S\\ \qp\le
T}}-b_{\pp}(\xx)\frac{\log(\qp)}{\qp}\right)\\&
=\res_{s=2}\left(\frac d{ds}(\log(L_2(\xx/k,s)))\right)
=\ord_{s=2}(L_2(\xx/k,s))
\end{align*}
    et
\begin{align*}&\lim_{T\to\infty}\frac
1T\left(\sum_{\substack{\pp\notin S\\ \qp\le
T}}-\tr(F_{\pp}\,|\,\ssc)\log(\qp)\right)\\ &=\res_{s=2}\left(
\frac d{ds}(\log(L(\ssc_{\ell}(1),s)))\right)
=\ord_{s=2}(L(\ssc_{\ell}(1),s)).\end{align*}
\newpage

En utilisant ces derni\`eres formules, la Proposition \ref{anal1a} et
la Conjecture \ref{conjT}, on peut
conclure
\begin{align*}&\lim_{T\to\infty}\frac
1T\left(\sum_{\substack{\pp\notin S\\ \qp\le
T}}-A_{\pp}^*(\xx)\log(\qp)\right)=-\ord_{s=2}(L_2(\xx/k,s))+\ord_{s=2
}(L(\ssc_{\ell}(1),s))\\
&=-\ord_{s=2}(N(\xx/k,s))=\rg\left(\frac{J_X(K)}{\tau B(k)}\right).\end{align*}
\end{proof}

\section{Exemples}

Les exemples d'application potentielle du r\'esultat sont tr\`es
nombreux \`a cause de la proposition
classique suivante.

\begin{proposition}\label{propU} Soit $\mathcal{X}$ une vari\'et\'e
lisse et projective de dimension $n$, d\'efinie sur un corps infini
$k$. Quitte \`a remplacer $\mathcal{X}$ par
son
\'eclatement en un nombre fini de points, il existe une fibration en
courbes de genre $g\geq 1$,
disons
$f:\mathcal{X}\rightarrow\ppb^{n-1}$, \'egalement d\'efinie sur $k$.
\end{proposition}

\begin{proof} Pour le cas des surfaces,  la construction d'un pinceau
de Lefschetz r\'epond \`a la question (voir par exemple \cite[Chapter
V, Proposition 3.1]{mi2}). Dans le
cas g\'en\'eral o\`u $\dim(\mathcal{X})=n$, prenons un plongement
$\mathcal{X}\subset\ppb^N$ et une
projection lin\'eaire
(\`a la Noether) induisant un morphisme fini
$\phi:\mathcal{X}\rightarrow\ppb^n$. Soit $p_0$ un point de
$\ppb^n$ situ\'e hors du diviseur de ramification de $\phi$; soit
$\pi:\tilde{\ppb}\rightarrow\ppb^n$ l'\'eclatement du point
$p_0$ et $\pi':\tilde{\ppb}\rightarrow\ppb^{n-1}$ l'application
associ\'ee (c'est un fibr\'e en droites
projectives). Si $\tilde{\mathcal{X}}$ est l'\'eclat\'e de
$\mathcal{X}$ au-dessus des points de
$\phi^{-1}(p_0)$, on obtient un
morphisme fini $\tilde{\phi}:\tilde{\mathcal{X}}\rightarrow
\tilde{\ppb}$ qui, compos\'e avec $\pi'$, fournit
la fibration en courbes cherch\'ee $f=\pi'\circ
\tilde{\phi}:\tilde{\mathcal{X}}\rightarrow
\ppb^{n-1}$.
\end{proof}

\begin{remark}
   En reprenant la construction pr\'ec\'edente, on voit qu'on a
\'eclat\'e $\mathcal{X}$ en $d$ points (avec
$d:=\deg_{\ppb^N}\mathcal{X}$). Si l'on consid\`ere $E$ l'un des
diviseurs exceptionnels sur $\tilde{\mathcal{X}}$, il est isomorphe
\`a $\ppb^{n-1}$ et le morphisme
$f_{|E}:E\rightarrow\ppb^{n-1}$ est de degr\'e 1 donc est un
isomorphisme. La construction fournit donc $d$ sections
de
$f:\tilde{\mathcal{X}}\rightarrow \ppb^{n-1}$ (d\'efinies sur une
   extension finie de $k$).
\end{remark}

Pour illustrer notre r\'esultat, il convient de choisir des exemples
de surfaces $\mathcal{X}$ pour lesquelles
la Conjecture de Tate (Conjecture \ref{conjT}) est d\'emontr\'ee.
Commen\c{c}ons par la proposition
suivante qui est bien
connue (par exemple l'argument est donn\'e dans \cite[Theorem
1.8]{rosi}, pour $\mathcal{X}=\ppb^2$).

\begin{proposition}\label{propTT} Soit $\mathcal{X},\mathcal{X}'$
deux surfaces lisses et projectives d\'efinies sur $k$.
   Si $\mathcal{X}$ et $\mathcal{X}'$ sont $k$-birationnellement
\'equivalente alors la conjecture de Tate est vraie pour
$\mathcal{X}/k$ si et seulement si elle est vraie
pour
$\mathcal{X}'/k$.
\end{proposition}

\begin{proof} En utilisant qu'une application birationnelle entre
surfaces lisses  est une composition de transformations monoidales et
de ses inverses \cite[Chapter V]{har}, on
se ram\`ene au cas o\`u on a un morphisme $\pi:\xx'\rightarrow\xx$
compos\'e  d'un nombre fini
d'\'eclate\-ments. Chacune
de ces transformations rajoute un diviseur exceptionnel $E$ au groupe
de N\'eron-Severi et au groupe de
cohomologie, c'est-\`a-dire
$\ns(\xx')\cong\pi^*(\ns((\xx))\oplus\mathbb{Z}[E]$ \cite[Proposition
II.3]{beau} et
$H^2_{\et}(\overline{\xx}',\mathbb{Q}_{\ell}))\cong
\pi^*(H^2_{\et}(\overline{\xx},\mathbb{Q}_{\ell}))\oplus\qq_{\ell}[E]$
\cite[p. 255]{rosi}. Globalement on
obtient donc, en notant $V$ le groupe des classes de diviseurs
engendr\'e par les diviseurs
exceptionnels,
$$\ns(\mathcal{X}')\cong\ns(\mathcal{X})\oplus V\text{ et }
H^2_{\text{\'et}}(\overline{\xx}',\qqell(1)))\cong
H^2_{\text{\'et}}(\overline{\xx},\qqell(1)))
\oplus
\left(V\otimes\qq_{\ell}(1)\right)$$ Le quotient
$L_2(\mathcal{X}'/k,s)/L_2(\mathcal{X}/k,s)=L(V,s-1)$ est donc
une fonction $L$ d'Artin en $s=2$ (cf. Proposition \ref{arbr}) d'ordre
\begin{equation}\label{tate01}\begin{aligned}
-\ord_{s=2}(L(V,s-1))&=\dim((V\otimes\qqell(1))^{\gk})\\
&=\dim(H^2_{\et}(\xbar',\qqell(1))^{\gk})-
\dim(H^2_{\et}(\xbar,\qqell(1))^{\gk}).\end{aligned}\end{equation}
D'autre part on a bien
$$\rg\left(\ns(\mathcal{X'}/k)\right)-\rg\left(\ns(\mathcal{X}/k)\right)
=\dim((V\otimes\qqell(1))^{\gk})$$ d'o\`u le r\'esultat.
\end{proof}

Dans tous les cas o\`u elle est d\'emontr\'ee,   la Conjecture
\ref{conjT} peut en fait se d\'ecompo\-ser  en
deux parties (cf. \cite{ram}, \cite{rosi}) dont la conjonction
entra\^{\i}ne visiblement le
r\'esultat cher\-ch\'e:
\begin{itemize}\item (\textbf{Conjecture de Tate I})
\begin{equation}\label{tate02}
\ns(\mathcal{X}/k)\otimes\qq_{\ell}\cong
H^2_{\text{\'et}}(\overline{\xx},\qq_{\ell}(1))^{\gk}\end{equation}
\item (\textbf{Conjecture de Tate II}) \emph{On a l'\'egalit\'e}
\begin{equation}\label{tate03} -\ord_{s=2}L_2(\mathcal{X}/k,s)=\dim
H^2_{\text{\'et}}(\overline{\xx},\qq_{\ell}(1))^{G_k}.\end{equation}
\end{itemize}

L'article de Rosen et Silverman \cite{rosi} d\'ecrit de nombreux exemples de surfaces (elliptiques) v\'erifiant la
Conjecture \ref{conjT}. Nous rassemblons dans le remarque suivant la plupart des r\'esultats connus (cf. \cite{ram},
\cite{rosi}) en les ordonnant selon la classification des surfaces, notamment  selon leur {\it dimension de Kodaira}
not\'ee $\kappa(\mathcal{X})$ (cf. \cite{bad} et \cite{beau}).

\begin{remark}\label{Tliste}
Soit $\mathcal{X}$ une surface lisse et projective d\'efinie sur un
corps de nombres $k$. La Conjecture de Tate I
est v\'erifi\'ee dans les cas suivants:
\begin{enumerate}
\item  La surface est rationnelle ou r\'egl\'ee
($\kappa(\mathcal{X})=-\infty$).

\item La surface $\mathcal{X}$ est une surface K3, une surface
d'Enriques, une surface ab\'elienne ou une
surface bi-elliptique($\kappa(\mathcal{X})=0$).

\item   La surface est le produit d'une courbe de genre 1 par une
courbe  de genre $\geq 2$ ($\kappa(\mathcal{X})=1$).

\item  La surface est un produit de courbes de genre $\geq 2$, une
surface modulaire de Hilbert ou une surface
de Fermat ($\kappa(\mathcal{X})=2$).
\end{enumerate}

La Conjecture de Tate II et donc la conjecture (\ref{conjT}) sont
\'egalement v\'erifi\'ees dans tous les cas
les suivants (a fortiori la conjecture de Nagao g\'en\'eralis\'ee est
v\'erifi\'ee pour toute
fibration d'une de ces
surfaces au-dessus d'une courbe (versions analytique et Taub\'erienne)):

\begin{enumerate}
\item La surface est rationnelle ou r\'egl\'ee ($\kappa(\mathcal{X})=-\infty$).

\item La surface est une surface K3 singuli\`ere ou de ``type CM",
une surface ab\'elienne CM ou facteur
d'une Jacobienne de courbe modulaire ($\kappa(\xx)=0$).

\item Une surface qui est produit d'une courbe modulaire de genre 1
par une courbe modulaire de genre $\ge
2$ ($\kappa(\xx)=1$).

\item  La surface est une surface modulaire de Hilbert ou une surface de Fermat ($\kappa(\mathcal{X})=2$).
\end{enumerate}

Dans le cas o\`u $\kappa(\xx)=-\infty$ les classes de diviseurs
engendrent $H^2_{\et}(\ov{\xx},\qq_{\ell}(1))$,
donc l'application des classes est d\'ej\`a un isomorphisme, en
particulier les Conjectures I et II
de Tate sont
satisfaites. Nous renvoyons \`a \cite[p. 237]{ram} pour l'affirmation
(et les r\'ef\'erences) que la Conjecture de
Tate I est connue pour les vari\'et\'es ab\'eliennes (Faltings), les
surfaces K3 (Deligne), les
surfaces modulaires de
Hilbert (Murty-Ramakrishnan) et les surfaces de Fermat (Shioda).

Les r\'esultats de Faltings entra\^{\i}nent \'egalement que la
Conjecture de Tate I est vraie pour un produit
de courbes (voir par exemple \cite[Theorem~5.2]{ta3}). Ajoutons que, si
$\mathcal{X}'\rightarrow\mathcal{X}$ est un
morphisme surjectif, alors la Conjecture de Tate I pour
$\mathcal{X}'$ entra\^{\i}ne relativement ais\'ement
la Conjecture de Tate I pour $\mathcal{X}$ (voir
\cite[Theorem~5.2]{ta3}). Or les surfaces
d'Enriques (resp. les surfaces
bi-elliptiques) poss\`edent un rev\^etement par une surface K3 (resp.
une surface ab\'elienne). La Conjecture de Tate
I est donc v\'erifi\'ee dans tous les exemples cit\'es.

Nous renvoyons \`a \cite[p. 242]{ram} pour l'affirmation (et les
r\'ef\'erences) que la Conjecture de Tate II
est connue pour les surfaces ab\'eliennes de type CM ou facteurs de
jacobiennes de courbes
modulaires (d\^u pour
l'essentiel \`a Shimura), les surfaces modulaires de Hilbert
(Murty-Ramakrishnan) et les surfaces de Fermat
(Shioda). Le cas des surfaces K3 singuli\`eres est d\^u \`a Shioda,
celui des surfaces K3 de ``type
CM" peut se d\'eduire du cas
des surfaces ab\'eliennes. Le cas de produit de courbes modulaires
peut \^etre trait\'e par la m\'ethode de
Rankin. Dans chaque cas on r\'eussit \`a exprimer
$L_2(\mathcal{X}/k,s)$ comme produit de fonctions
$L$ associ\'ees \`a une
forme modulaire ou \`a un caract\`ere de Hecke et on sait donc que
$L_2$ ne s'annule pas sur
$\Re(s)=2$.
\end{remark}

\begin{remark}
\begin{enumerate}
\item On appelle ici surface K3 de ``type CM" celles pour les\-quel\-les le motif est de type CM, i.e.,
   sa fonction $L$
est produit de fonctions $L$ associ\'ees \`a des {\it
Gr\"ossencharacter} de Hecke. Pour une vari\'et\'e ab\'elienne
de type CM, Pohlman montre dans \cite{pol} que la conjecture de Tate
est \'equivalente \`a la
conjecture de Hodge pour ces
vari\'et\'es. De plus la conjecture de Hodge est vraie  en
codimension 1, d'apr\`es un r\'esultat de Lefschetz.
En particulier, elle est vraie pour les surfaces ab\'eliennes de type
CM. L'argument pour les
surfaces K3 de type CM est
similaire \`a \cite{pol}.

\item Les surfaces de Fermat sont les hypersurfaces de $\ppb^3$ donn\'ees par une \'equa\-tion
$a_0x_0^d+a_1x_1^d+a_2x_2^d+a_3x_3^d=0$ (si $d\leq 3$ c'est une surface rationnelle, si $d=4$ c'est une surface K3 de
``type CM" et si $d\geq 5$, c'est une surface de type g\'en\'eral avec $\kappa(\mathcal{X})=2$).
\end{enumerate}

\end{remark}

\begin{example}\label{exsh} Donnons un exemple concret en genre $g$
quelconque, essentiellement tir\'e de \cite{sh1}. Soit
$P(x)=x^{2g+1}+p_{2g}x^{2g}+\dots+p_0$ un polyn\^ome
s\'eparable de $k[x]$, on d\'efinit alors la courbe $X$ sur $k(t)$
comme la courbe lisse projective
dont un ouvert affine est
donn\'e par
$$y^2=x^{2g+1}+p_{2g}x^{2g}+\dots+p_0+t^2.$$ Si
$f:\mathcal{X}\rightarrow\ppb^1$ est une fibration d'une surface
projective lisse d\'efinie sur $k$ dont la
fibre g\'en\'erique est $X$, on voit que $\mathcal{X}$ est
$k$-rationnelle (en effet la surface est
birationnelle \`a la
surface d'\'equation $uv=P(x)$). Comme $\pic^0(\mathcal{X})=0$, on a
aussi $B=0$. Soit $O$ le point \`a l'infini de
$X$ et $P_i=(e_i,t)$ (respectivement $P'_i=(e_i,-t)$) les points de
$X(\bar{k}(t))$ correspondant
aux z\'eros de
l'\'equation $P(x)=0$. On voit tout de suite que, si on plonge $X$
dans $J_X$ en prenant $O$ comme point base, on a
les relations $P_i+P'_i=0$ et $P_1+\dots+P_{2g+1}=0$. On sait par
Shioda \cite{sh1} que
$\rg(J_X(\bar{k}(t)))=2g$ et que
les points $P_1,\dots,P_{2g}$ forment un syst\`eme g\'en\'erateur.
Supposons de plus $e_i\in k$, ou encore
$P_i\in X(k(t))$ alors on a aussi $\rg(J_X(k(t)))=2g$ et donc, dans
l'exemple pr\'esent
$A_{\pp}^*(\xx)=A_{\pp}(\xx)$ et on a


$$
\lim_{T\to\infty}\frac 1T\left(\sum_{\substack{\pp\notin S\\ \qp\le
T}}-A_{\pp}(\xx)\log(\qp)\right)=2g.$$ Si l'on
ne suppose plus que $P$ a toutes ses racines dans $k$ on doit
remplacer $2g$ par le nombre de
facteurs  moins un dans la
d\'ecomposition en facteurs irr\'eductibles de $P$ dans $k[x]$.
\end{example}

\begin{example} Donnons un autre exemple avec cette fois une
$\qq(t)/\qq$-trace non nulle, ce qui rend donc indispensable de
consid\'erer la trace moyenne r\'eduite des
Frobenius.
   Soit $X$ la courbe lisse
projective d\'efinie sur $\qq(t)$ dont un ouvert affine est donn\'e ainsi~:
$$U:=\{(x,y,z)\in\aab^3\;|\;y^2-(x-a)(x-b)(x-c)=z^2-x(x-1)(x-t)=0\},$$
o\`u $a,b,c\in\qq^*-\{1\}$ sont distincts. Soit $V$ la surface affine
lisse d\'efinie sur $\qq$
ainsi~:
$$V:=\{(x,y,z,t)\in\aab^4\;|\;y^2-(x-a)(x-b)(x-c)=z^2-x(x-1)(x-t)=0\quad\hbox{et}\quad t\not= 0\}$$
et munie de la fibration $f(x,y,z,t)=t$ de $V$ vers
$\aab^1\setminus\{0\}$. On peut construire une surface
lisse projective $\xx/\qq$ munie d'une fibration not\'ee encore $f$
de $\xx$ vers $\ppb^1$ telle
que sa restriction \`a $V$
soit la fibration $f:V\rightarrow \aab^1\setminus\{0\}$. La fibre
g\'en\'erique de $f:\xx\rightarrow\ppb^1$
est naturellement $X/\qq(t)$. Introduisons  les courbes lisses
projectives  d\'efinies par
l'\'equation d'un ouvert affine
\begin{multline*}
E_1/\qq: \, y^2=(x-a)(x-b)(x-c),\; E_2/\qq(t): \,
z^2=x(x-1)(x-t)\quad\hbox{et}\\
C_1/\qq(t): \, w^2=x(x-1)(x-t)(x-a)(x-b)(x-c).
\end{multline*}

   Le genre de $E_1$ et $E_2$ est 1, le genre
de $C_1$ est $2$ et le genre de $X$ est 4. L'application
$\pi:V\rightarrow E_1$ donn\'ee par $\pi(x,y,z,t)=(x,y)$
a pour fibres des coniques. Ainsi la surface $\xx$ est birationnelle
\`a une surface r\'egl\'ee au
dessus de $E_1$ et la
conjecture de Tate est donc v\'erifi\'ee pour cette surface; de plus
ceci montre aussi que $E_1$ est  la
vari\'et\'e d'Albanese de $\xx$ et par cons\'equent \'egalement sa
$\qq(t)/\qq$-trace car
$E_1=E_1^{\vee}$. On peut montrer que
$J_X$ est isog\`ene \`a $E_1\times E_2\times J_{C_1}$. On a donc
$$\rg(J_X(\qq(t))/E_1(\qq))=\rg(E_2(\qq(t)))+\rg(J_{C_1}(\qq(t))).$$
Il est bien connu que $\rg(E_2(\qq(t)))=0$. En utilisant une
2-descente, on peut montrer que
$\rg(J_{C_1}(\qq(t)))=1$ (avec comme g\'en\'erateur d'un sous-groupe
d'indice fini la diff\'erence
des deux points \`a l'infini de $C_1$), pour
un choix assez g\'en\'eral de $a,b,c\in\qq^*$. Si l'on pose comme
pr\'ec\'edemment
$$A_{\pp}^*(\xx):=\frac
1{\qp}\left(\sum_{c\in\ppb^1(\kappa_{\pp})}a_{\pp}(\xx_c)\right)
-a_{\pp}(E_1),$$ on obtient donc en appliquant
le Th\'eor\`eme \ref{thmA}:

$$
\lim_{T\to\infty}\frac 1T\left(\sum_{\substack{\pp\notin S\\ \qp\le T}}
-A_{\pp}^*(\xx)\log(\qp)\right)=1.$$
\end{example}

Nous terminerons en notant deux corollaires du Th\'eor\`eme~\ref{thmA}.

\begin{corollary}\label{corTa-Na} La conjecture de Tate pour les
surfaces est \'equivalente \`a la conjecture analytique de Nagao
(g\'en\'eralis\'ee).
\end{corollary}

\begin{proof} La conjecture de Tate entra\^{\i}ne la conjecture
analytique de Nagao d'apr\`es le Th\'eor\`eme \ref{thmA}. Inversement
si $\mathcal{X}/k$ est une surface, elle
est birationnelle \`a $\mathcal{X}'/k$ poss\'edant une fibration en
courbes sur $\ppb^1$, donc la
conjecture de Nagao
analytique entra\^{\i}ne la conjecture de Tate pour $\mathcal{X}'/k$
(cf. Proposition \ref{anal1a}) et donc
pour
$\mathcal{X}/k$ d'apr\`es la Proposition \ref{propTT}.
\end{proof}

\begin{corollary}\label{corNa1-2} La conjecture  analytique de Nagao
(g\'en\'eralis\'ee) pour les  fibrations au dessus de $\ppb^1$
entra\^{\i}ne la conjecture analytique de
Nagao (g\'en\'eralis\'ee) pour les  fibrations au dessus d'une courbe
quelconque $C$.
\end{corollary}

\begin{proof} D'apr\`es l'argument pr\'ec\'edent, la conjecture pour
les fibrations sur $\ppb^1$ suffit \`a entra\^{\i}ner la conjecture
de Tate, et cette derni\`ere entra\^{\i}ne
la conjecture analytique de Nagao.
\end{proof}

\end{document}